\documentclass{svjour3}

\usepackage[utf8]{inputenc}

\usepackage{amsmath}
\usepackage{amssymb}
\usepackage{mathtools}
\usepackage{xcolor}
\usepackage{listings}
\usepackage{pgfplots}\pgfplotsset{compat=1.16}

\usepackage{hyperref}
\usepackage{hyperxmp}
\colorlet{internallinkcolour}{green!50!black}
\colorlet{externallinkcolour}{red!50!black}
\hypersetup
{
  pdfauthor          = {Matthias Kirchhart; Donat Weniger},
  pdfsubject         = {Preprint},
  pdftitle           = {Analytic Integration of the Newton Potential over Cuboids
and an Application to Fast Multipole Methods},
  pdfkeywords        = {Newton potential, analytic integration, Poisson equation, C++},
  pdfencoding        = {unicode},
  pdflang            = {en-GB},
  pdfdisplaydoctitle = {true},
  unicode            = {true},
  colorlinks         = {true}, 
  allcolors          = internallinkcolour,
  urlcolor           = externallinkcolour,
  frenchlinks        = false,
  pdfborder          = {0 0 0},
  naturalnames       = false,
  hypertexnames      = false,
  breaklinks         = true
}

\usepackage
[
	backend=biber,
	language=british,
	sorting=none,
    sortcites=true,
	sortlocale=en-GB,
    giveninits=true,
    maxnames=4,
    isbn=false,
    doi=false,
    autocite=superscript
]{biblatex}
\DeclareCiteCommand{\supercite}[\mkbibsuperscript]{%
\iffieldundef{prenote}{}{\BibliographyWarning{Ignoring prenote argument}}%
\iffieldundef{postnote}{}{}}%
{\bibopenbracket\usebibmacro{citeindex}\usebibmacro{cite}\usebibmacro{postnote}%
\bibclosebracket}{\supercitedelim}{}

\DeclareMathOperator\atanh{artanh}
\DeclareMathOperator\atan{arctan}
\DeclareMathOperator\supp{supp}

\newcommand{\bigO}[1]{\mathcal{O}\left({#1}\right)}
\newcommand{\spdim}{\mathrm{D}}   
\newcommand{\reals}{\mathbb{R}}
\newcommand{\wholespace}{\reals^\spdim}

\newcommand{\dd}{{\mathrm d}}

\newcommand{\vv}[1]{\mathbf{#1}}        
\newcommand{\gv}[1]{\boldsymbol{#1}}    

\begin{document}
\titlerunning{Analytic Integration of the Newton Potential}
\title{Analytic Integration of the Newton Potential over Cuboids and an
Application to Fast Multipole Methods}
\author{Matthias Kirchhart \and Donat Weniger}
\institute{Applied and Computational Mathematics, RWTH Aachen University, Schinkelstra{\ss}e 2,
52062~Aachen, Germany\\\email{kirchhart@acom.rwth-aachen.de \and
donat.weniger@rwth-aachen.de}}
\maketitle

\begin{abstract}
We present simplified formulæ for the analytic integration of the Newton
potential of polynomials over boxes in two- and three-dimensional space. These
are implemented in an easy-to-use C++ library that allows computations in
arbitrary precision arithmetic which is also documented here. We describe how
these results can be combined with fast multipole methods for general, non-%
polynomial data.
\end{abstract}


\section{Introduction}
Consider the Poisson equation on the whole space:
\begin{equation}\label{eqn:poisson-equation}
-\Delta u = f\quad\text{on }\wholespace
\end{equation}
where $f$ is a compactly supported function and $\mathrm{D}\in\lbrace 2,3\rbrace$.
It is well-known that a solution to this problem is given by the Newton potential.
To this end, let us define the fundamental solution $G$:
\begin{equation}
G: \wholespace\to\reals,\quad \vv{x}\mapsto
\left\lbrace
\begin{aligned}
\frac{1}{2\pi}&\ln\left(\frac{1}{|\vv{x}|}\right) & \spdim = 2, \\
\frac{1}{4\pi}&\frac{1}{|\vv{x}|}                 & \spdim = 3.
\end{aligned}
\right.
\end{equation}
Then a solution to~\eqref{eqn:poisson-equation}
 -- the Newton potential of $f$ -- 
is given by $u=G\star f$, that is:
\begin{equation}
\label{eqn:convolution}
u(\vv{x}) = \int_{\wholespace}G(\vv{x}-\vv{y})f(\vv{y})\,{\mathrm d}{\vv{y}}.
\end{equation}

Note that whenever $\vv{x}\in\supp f$ this integral is singular and standard
quadrature methods become ineffective. Now suppose $f$ would be given as a
piece-wise polynomial on axis-aligned boxes. In this case the difficulty due to
the singularity can be avoided by computing this integral analytically. If
moreover the orthogonal projection $Pu$ of $u=G\star f$ onto the same space is
sought, another integration needs to be carried out. In three dimensions, for
example, this leads to integrals of the form:
\begin{equation}
  \label{eqn:newton-potential}
    \iiint\displaylimits_{Q'} \iiint\displaylimits_{Q}
    \frac{x_1^{\mu_1} x_2^{\mu_2} x_3^{\mu_3} y_1^{\lambda_1} y_2^{\lambda_2} y_3^{\lambda_3}}
         {\sqrt{\left(x_1 - y_1\right)^2 + \left(x_2 - y_2\right)^2 + \left(x_3 - y_3\right)^2}}
    \,{\dd y_1}{\dd y_2}{\dd y_3}{\dd x_1}{\dd x_2}{\dd x_3}
\end{equation}
with arbitrary monomials in the enumerator, so $\gv{\mu}=(\mu_1,\mu_2,\mu_3)$,
$\gv{\lambda}=(\lambda_1,\lambda_2,\lambda_3)\in\mathbb{N}_0^3$ are multi-indices.
The integration domains are cuboids, $Q = \prod\limits_{i=1}^3\left[a_i,b_i\right]$
and $Q' = \prod\limits_{i=1}^3\left[a_i',b_i'\right]$. Again, whenever $Q\cap Q'
\neq\emptyset$ the singularity prohibits the use of standard quadrature
formulæ.

For the three-dimensional case analytical integration formulæ for~%
\eqref{eqn:newton-potential} have been presented by Hackbusch.\autocite{hackbusch2002,hackbusch2002b}
His formulæ are, however, very complicated and difficult to implement because
they distinguish between different cases and involve several intermediate
functions. He furthermore describes difficulties from numerical cancellation
in certain cases. In this work we present formulæ in terms of the $\atanh$
function, which are not only significantly simpler, but also slightly less
prone to cancellation errors. In addition to Hackbusch's work, we also present
the corresponding two-dimensional formulæ.

Given the availability of modern computer algebra systems one might be inclined
to ask whether manually tabulated integrals as the ones in subsection~\ref{ch:basic_integrals_3d}
are still adequate. In our experiments however we realised that these systems
give formulæ that are simply over-complicated and usually fail after the second
or third integration. Obtaining the simplicity and symmetry of the formulæ below
is a tedious and error-prone process from which we would like to spare others.

To avoid the effect of numerical cancellation, the formulæ were implemented in a
simple C++ library, which performs the evaluation in arbitrary precision. We hope
that this software is useful for other researchers and therefore release it as
free software under the terms of the GNU Lesser General Public License, version
three or later.

In numerical software, these analytic integration formulæ for the singular
near-field integrals can be combined with fast multipole methods or tree codes
for the far-field computation. In case of uniform Cartesian grids, the
integration results can be pre-computed once in high-precision for a single
reference cell interacting with its neighbours. Because convolution commutes
with translation, the results can then be reused for all other near-field
interactions. This was for example used recently in a vortex method,\autocite{kirchhart2019b}
and we give a description of such a scheme at the end.

Finally, let us remark that there also exist specialised, numerical quadrature
schemes for this problem, for example the well-known Duffy trick\autocite{duffy1982}
or hierarchical quadrature.\autocite{hackbusch2005,hackbusch2008} These techniques
have the advantage of also generalising to other integration kernels $G$, but are
usually not exact. The software presented here can be used to obtain arbitrarily
accurate reference data for testing these and other similar schemes.

\section{The Three-dimensional Case}

\subsection{Shape of Antiderivatives}
Similar to Hackbusch, we make use of the following abbreviations:\autocite{hackbusch2002}
\begin{gather}
    X\coloneqq x_1 - y_1,\quad
    Y\coloneqq x_2 - y_2,\quad
    Z\coloneqq x_3 - y_3, \\
    R\coloneqq\sqrt{X^2+Y^2+Z^2}.
\end{gather}
The integral \eqref{eqn:newton-potential} can be evaluated by computing the
six-fold antiderivative of the integrand with respect to the coordinate directions
$x_1$, $x_2$, $x_3$, and $y_1$, $y_2$, $y_3$. Once this is done, the evaluation of the
integral~\eqref{eqn:newton-potential} amounts to evaluating the resulting function at
the integration boundaries. These anti-derivatives are always of the form:
\begin{multline}\label{eqn:antiderivative3d}
    F(x_1,x_2,x_3,y_1,y_2,y_3)\coloneqq
      P_1 \frac{1}{R}
    + P_2 R \\
    + P_3 \atanh{\left(\frac{X}{R}\right)} 
    + P_4 \atanh{\left(\frac{Y}{R}\right)}
    + P_5 \atanh{\left(\frac{Z}{R}\right)} \\ 
    + P_6 \atan{\left(\frac{X}{R}\frac{Y}{Z}\right)}
    + P_7 \atan{\left(\frac{X}{R}\frac{Z}{Y}\right)}
    + P_8 \atan{\left(\frac{Y}{R}\frac{Z}{X}\right)},
\end{multline}
where $P_i$, $i=1,\ldots,8$ are certain polynomials with rational coefficients
$c_{i,\gv{\lambda},\gv{\mu}}\in\mathbb{Q}$. It turns out that integration is
closed over functions $F$ of this type: integrating $F$ with respect to one of
the coordinate directions will always yield a new function $\widetilde{F}$ of
the same shape. Interestingly, unlike in Hackbusch's formulation, there is no
stand-alone polynomial term.

The $P_1$-term always vanishes after a single integration. After integrating over
$y_1$, $y_2$, and $y_3$, or $x_1$, $x_2$, and $x_3$, the polynomials $P_3, P_4,
\ldots, P_8$ will have a multiple zero whenever the corresponding $\atanh$ or
$\atan$ terms are singular. With the singularities cancelled out in this way,
the overall resulting antiderivative $F$ for the computation of
$\eqref{eqn:newton-potential}$ will always be a globally continuous function.

\subsection{Basic Integrals}
\label{ch:basic_integrals_3d}
We first consider the case where all but one of the $P_i,i=1,\ldots,8$ are zero,
with the remaining polynomial $P_j\equiv 1$. The results for integration with
respect to $x_1$ are as follows:
\begin{align}
  \begin{split}
  \label{eqn:antiderivative-1}
    \int \frac{1}{R}  \, \dd x_1 = \atanh{\left(\frac{X}{R}\right)},
  \end{split} \\
  \begin{split}
      \int R  \, \dd x_1 = \frac{1}{2}X R + \frac{1}{2} \left(Y^2 + Z^2\right) \atanh{\left(\frac{X}{R}\right)},
  \end{split}
\end{align}
\begin{align}
  \begin{split}
    \int \atanh{\left(\frac{X}{R}\right)}  \, \dd x_1 = X \atanh{\left(\frac{X}{R}\right)} - R,
  \end{split} \\
  \begin{split}
    \label{eqn:antiderivative-4}
    \int \atanh{\left(\frac{Y}{R}\right)}  \, \dd x_1 = X \atanh{\left(\frac{Y}{R}\right)}
      + Y \atanh{\left(\frac{X}{R}\right)}
      - Z \atan{\left(\frac{Y}{R}\frac{X}{Z}\right)},
  \end{split} \\
  \begin{split}
    \label{eqn:antiderivative-5}
    \int \atanh{\left(\frac{Z}{R}\right)}  \, \dd x_1 = X \atanh{\left(\frac{Z}{R}\right)}
      + Z \atanh{\left(\frac{X}{R}\right)}
      - Y \atan{\left(\frac{Z}{R}\frac{X}{Y}\right)},
  \end{split}
\end{align}
\begin{align}
  \begin{split}
    \int \atan{\left(\frac{X}{R}\frac{Y}{Z}\right)}  \, \dd x_1 = X \atan{\left(\frac{X}{R}\frac{Y}{Z}\right)}
      + Z \atanh{\left(\frac{Y}{R}\right)},
  \end{split} \\
  \begin{split}
    \int \atan{\left(\frac{X}{R}\frac{Z}{Y}\right)}  \, \dd x_1 = X \atan{\left(\frac{X}{R}\frac{Z}{Y}\right)}
      + Y \atanh{\left(\frac{Z}{R}\right)},
  \end{split} \\
  \begin{split}
    \label{eqn:antiderivative-8}
    \int \atan{\left(\frac{Y}{R}\frac{Z}{X}\right)}  \, \dd x_1 = X \atan{\left(\frac{Y}{R}\frac{Z}{X}\right)}
      - Y \atanh{\left(\frac{Z}{R}\right)}
      - Z \atanh{\left(\frac{Y}{R}\right)}.
  \end{split}
\end{align}
The integration for the variables $x_2$ and $x_3$ are equivalent modulo slight
rearrangements of $X$, $Y$ and $Z$. For the antiderivatives with respect to $x_2$,
one has to switch $X$ and $Y$ in formulæ \eqref{eqn:antiderivative-1} to
\eqref{eqn:antiderivative-8}, and adjust the differential to $\dd x_2$. Applying
this rule to~\eqref{eqn:antiderivative-4} for example yields:
\begin{equation}
    \int \atanh{\left(\frac{X}{R}\right)}  \, \dd x_2 = Y \atanh{\left(\frac{X}{R}\right)}
      + X \atanh{\left(\frac{Y}{R}\right)}
      - Z \atan{\left(\frac{X}{R}\frac{Y}{Z}\right)}.
\end{equation}
Analogously, for the $x_3$ integration $X$ and $Z$ are switched in formulæ
\eqref{eqn:antiderivative-1} to \eqref{eqn:antiderivative-8} and the differential
is changed to $\dd x_3$.

The integrations for $y_1$, $y_2$ and $y_3$ are also easily extracted since the
$\atanh{}$ and $\atan{}$ functions are odd. The antiderivatives with respect to
$y_i$ are therefore the same as the ones with respect to $x_i$ multiplied by $-1$,
$i \in\lbrace 1, 2, 3\rbrace$. Using this for the $y_1$ integration of the term in,
e.\,g., formula~\eqref{eqn:antiderivative-5}, we get
\begin{equation}
    \int \atanh{\left(\frac{Z}{R}\right)}  \, \dd y_1 = -X \atanh{\left(\frac{Z}{R}\right)}
      - Z \atanh{\left(\frac{X}{R}\right)}
      + Y \atan{\left(\frac{Z}{R}\frac{X}{Y}\right)}.
\end{equation}

\subsection{General Polynomials}
We now consider the case where the polynomials $P_i$, $i=1,\ldots,8$ are non-trivial.
Let us for example take the term $P_1/R$ and look at the integration with
respect to $x_1$. Writing $P_1=\sum_{\gv{\lambda},\gv{\mu}}c_{1,\gv{\lambda},\gv{\mu}}
\vv{x}^{\gv{\lambda}}\vv{y}^{\gv{\mu}}$ one obtains:
\begin{equation}
\int\frac{P_1}{R}\,{\mathrm d}x_1 =
\sum_{\gv{\lambda},\gv{\mu}} c_{1,\gv{\lambda},\gv{\mu}}
x_2^{\lambda_2}x_3^{\lambda_3}y_1^{\mu_1}y_2^{\mu_2}y_3^{\mu_3}
\int\frac{x_1^{\lambda_1}}{R}\,{\mathrm d}x_1.
\end{equation}
It thus suffices to consider monomials of the form $x_1^{\lambda_1}$. This task can
be solved using integration by parts.

For example, for $\lambda_1=1$ one obtains using the results of the previous section:
\begin{multline}
\int\frac{x_1}{R}\,{\mathrm d}x_1 = x_1\int\frac{\dd x_1}{R} -
\iint\frac{1}{R}\,\dd x_1 \dd x_1 \\
= x_1\atanh\left(\frac{X}{R}\right) - \int\atanh\left(\frac{X}{R}\right)\,\dd x_1 \\
= R + y_1\atanh\left(\frac{X}{R}\right).
\end{multline}
For the other functions $P_i$ and coordinate directions one can proceed in an
analogous fashion.

For general $\lambda_1\in\mathbb{N}$, integration by parts yields mutually nested
recursion relations. They are lengthy, but can be easily implemented on a computer.
The complete set of recurrence relations is available as supplementary material to
this article and has been implemented in the software library described below.

\subsection{A Remark on the Area Hyperbolic Tangent}
Formula \eqref{eqn:antiderivative-1} can, up to the integration constant,
be rewritten in terms of the natural logarithm as:
\begin{equation}
\int\frac{1}{R}\,\dd x_1 = \ln(X+R).
\end{equation}
This formula is problematic when implemented in floating point arithmetic.
When evaluating the \emph{sum} $X+R$ numerical cancellation occurs whenever
$|X|\approx R$ and $X<0$. The resulting error then gets further amplified
by the logarithm's singularity at zero.

The formulation in terms of the area hyperbolic tangent avoids this because
only the \emph{ratio} $X/R$ needs to be evaluated. The remaining singularity
when $|X|=R$ can be avoided in a numerical implementation by a simple check
for equality or by instead implementing:
\begin{equation}
\int\frac{1}{R}\,\dd x_1 \approx \atanh\left(\frac{X}{R}(1-c\,\mathrm{eps})\right),
\end{equation}
where $c>1$ is a small constant and ``$\mathrm{eps}$'' is the machine
epsilon.

\section{The Two-Dimensional Case}
The characteristics of the integration in three dimensions are reflected in the
two-dimensional case. We reuse the abbreviations
\begin{gather}
    X\coloneqq x_1 - y_1,\quad
    Y\coloneqq x_2 - y_2,
\end{gather}
and note that because 
\begin{math}
\frac{1}{2\pi}\ln\left(\frac{1}{|\vv{x}|}\right) = -\frac{1}{4\pi}\ln(|\vv{x}|^2)
\end{math}
it is sufficient to consider integrals of the form:
\begin{equation}
\int P_1 \ln{\left(X^2+Y^2\right)}\,{\mathrm d}\xi\qquad \xi\in\lbrace x_1,x_2,y_1,y_2\rbrace.
\end{equation}
The resulting antiderivatives will always be of the form:
\begin{equation}
    H(x_1,x_2,y_1,y_2)\coloneqq
      P_1 \ln{\left(X^2+Y^2\right)}
    + P_2 \atan{\left(\frac{X}{Y}\right)}
    + P_3 \atan{\left(\frac{Y}{X}\right)}
    + P_4.
\end{equation}
The $P_i, i=1,...,4$, are again polynomials with rational coefficients,
integration is closed over functions $H$ of this shape, and zeros in the
polynomials $P_1, P_2$ and $P_3$ cancel out the singularities in the
corresponding $\ln$ and $\atan$ terms. The corresponding basic integrals
with respect to $x_1$ read:
\begin{align}
  \begin{split}
    \int \ln{\left(X^2+Y^2\right)} \, \dd x_1 = X\ln{\left(X^2+Y^2\right)}-2X-2Y\atan{\left(\frac{Y}{X}\right)},
  \end{split} \\
  \begin{split}
    \int \atan{\left(\frac{X}{Y}\right)} \, \dd x_1 = X\atan{\left(\frac{X}{Y}\right)} - \frac{Y}{2}\ln{\left(X^2+Y^2\right)},
  \end{split} \\
  \begin{split}
    \int \atan{\left(\frac{Y}{X}\right)} \, \dd x_1 = X\atan{\left(\frac{Y}{X}\right)} + \frac{Y}{2}\ln{\left(X^2+Y^2\right)}.
  \end{split}
\end{align}

The integration for the remaining variables can be derived easily and is
available with the complete set of nested recursion formulæ as supplementary
material in addition to the implementation.

\section{Brief Overview of the C++ Library}
\lstset{basicstyle={\small\ttfamily},language=[11]C++,}
For brevity we will describe the structure of the library for the
three-dimensional case ($\spdim=3$) only. The two-dimensional implementation
is completely analogous.

\subsection{Data Structures}

\paragraph{Multi-indeces.} These are represented in an obvious way:
\begin{lstlisting}
template <size_t dim> using multi_index = std::array<size_t,dim>;
\end{lstlisting}
For $\spdim=3$ we choose \lstinline{dim==6}. A single instance
\lstinline{idx} represents a pair $(\gv{\lambda},\gv{\mu})\in
\mathbb{N}_0^6$. The first three entries \lstinline{idx[0]},
\lstinline{idx[1]}, \lstinline{idx[2]} correspond to $\lambda_1,
\lambda_2,\lambda_3$, and \lstinline{idx[3]},
\lstinline{idx[4]}, \lstinline{idx[5]} correspond to $\mu_1,\mu_2,
\mu_3$.

\paragraph{Polynomials.} These are stored as an associative
array of multi-indeces $(\gv{\lambda},\gv{\mu})$ and corresponding
coefficients $c_{\gv{\lambda},\gv{\mu}}\in\mathbb{Q}$. In the
software library, the user may specify the data type that is
used to represent the coefficients as a template parameter:
\begin{lstlisting}
template <typename number_type, size_t dim>
class polynomial;
\end{lstlisting}
One may choose \lstinline{double} for \lstinline{number_type}, however
this would come with a loss of precision. For exact symbolic computations
it is advisable to use rational arithmetic types such as
\lstinline{boost::multiprecision::mpq_rational} from Boost
and the GNU Multiprecision Library.\footnote{\url{https://www.boost.org},
\url{https://gmplib.org}}

The \lstinline{polynomial} class allows instances to be used in an
intuitive way using standard mathematical notation. The polynomial
$P = 16\,x_1y_2^3 + \frac{2}{3}\,x_1x_2^2$ may for example be created as
follows.
\begin{lstlisting}
using idx6     = multi_index<6>;
using rational = boost::multiprecision::mpq_rational;
using poly6    = polynomial<rational,6>;

// Short-hands for convenience.
const poly6 x1 = poly6( rational(1), idx6 { 1, 0, 0, 0, 0, 0 } );
const poly6 x2 = poly6( rational(1), idx6 { 0, 1, 0, 0, 0, 0 } );
const poly6 x3 = poly6( rational(1), idx6 { 0, 0, 1, 0, 0, 0 } );
const poly6 y1 = poly6( rational(1), idx6 { 0, 0, 0, 1, 0, 0 } );
const poly6 y2 = poly6( rational(1), idx6 { 0, 0, 0, 0, 1, 0 } );
const poly6 y3 = poly6( rational(1), idx6 { 0, 0, 0, 0, 0, 1 } );

poly6 P = 16*x1*y2*y2*y2 + rational(2,3)*x1*x2*x2;
\end{lstlisting}

\paragraph{Antiderivatives.} Functions of the form~\eqref{eqn:antiderivative3d}
are represented as objects of type \lstinline{antiderivative}, the
two-dimensional analogue is called \lstinline{antiderivative2d}.
\begin{lstlisting}
template <typename number_type>
class antiderivative
{
public:
    // ...
    using poly6  = polynomial<number_type,6>;

    poly6 PR, PRinv,
          Partanh_XR, Partanh_YR, Partanh_ZR,
          Parctan_RZ, Parctan_RY, Parctan_RX;
};
\end{lstlisting}
The \lstinline{poly6} members directly correspond to the
$P_i$ in equation \eqref{eqn:antiderivative3d}.
 
\subsection{Integration and Evaluation Routines}
Once an \lstinline{antiderivative} object has been created,
it may be integrated \emph{arbitrarily often} with respect to the
different coordinates. For example, continuing with the example
$P = 16\,x_1y_2^3 + \frac{2}{3}\,x_1x_2^2$, the following
code-snippet computes the indefinite integral:
\begin{equation}
\iiint \iiint
\frac{16\,x_1y_2^3 + \frac{2}{3}\,x_1x_2^2}
   {\sqrt{\left(x_1 - y_1\right)^2 + \left(x_2 - y_2\right)^2 + \left(x_3 - y_3\right)^2}}
   \,{\dd y_1}{\dd y_2}{\dd y_3}{\dd x_1}{\dd x_2}{\dd x_3}.
\end{equation}

\begin{lstlisting}
antiderivative<rational> F; // Zero by default.
F.PRinv = P;                // P as created above.

F = y1_integrate(F);
F = y2_integrate(F);
F = y3_integrate(F);
F = x1_integrate(F);
F = x2_integrate(F);
F = x3_integrate(F);
\end{lstlisting}
The order of the integrations is arbitrary, as in the exact mathematical
expression.

Once integration is complete, the resulting function \lstinline{F} may be
evaluated at arbitrary locations. This can be done using the member templates
\begin{lstlisting}
template<typename number_type>
class antiderivative
{
public:
    // ...
    template <typename float_type> 
    float_type eval( const std::array<number_type,6> &pos ) const;

    template <typename float_type, typename output_iterator>
    void write_summands( const std::array<number_type,6> &pos,
                         output_iterator out ) const;
};
\end{lstlisting}
Passing $(\vv{x},\vv{y})\in\mathbb{Q}^6$ will then result in an \emph{exact}
evaluation of the polynomials $P_i$ if a rational arithmetic type is used
for \lstinline{number_type}. Afterwards the result is converted and rounded
to the closest \lstinline{float_type}. This is necessary because the $\atanh$,
$\atan$, and $R$ terms are usually irrational numbers. The \lstinline{eval}
function then directly computes the result, whereas the \lstinline{write_summands}
method writes the resulting eight summands to the given output iterator. This
can be useful for estimating the condition number of the sum.

In continuing the example, to evaluate the resulting antiderivative at
$(\vv{x},\vv{y})=(1.234,1,2,3,4,5)$ using an accuracy of
$100$ decimal digits for \lstinline{float_type} one might use:
\begin{lstlisting}
using float100 = boost::multiprecision::mpf_float_100;
std::cout << std::setprecision(50) 
          << F.eval<float100>( { rational(1234,1000), 1, 2, 3, 4, 5 } );
\end{lstlisting}
which outputs the first 50 digits:
\begin{equation}
32801.793402730158138263112227926388836137928920324.
\end{equation}

\subsection{Effects of Numerical Cancellation}
In his work Hackbusch considers the following example:
\autocite[Section~6.10]{hackbusch2002b}
\begin{equation}
I \coloneqq
\smashoperator{\int_{0}^{100}}
\smashoperator{\int_{0}^{1}}  
\smashoperator{\int_{0}^{1}} 
\smashoperator{\int_{0}^{1}}
\smashoperator{\int_{0}^{100}}
\smashoperator{\int_{0}^{1}}
\frac{\dd y_1\dd y_2\dd y_3\ \dd x_1\dd x_2\dd x_3}{\sqrt{(x_1-y_1)^2+(x_2-y_2)^2+(x_3-y_3)^2}},
\end{equation}
and mentions ill-conditioning of the resulting sum. The condition
$\kappa$ of a sum $\sum_{k=1}^{N}a_k$ is defined as
\begin{equation}
\kappa\coloneqq\frac{\sum_{k=1}^{N}|a_k|}{\left|\sum_{k=1}^{N}a_k\right|}.
\end{equation}

As a rule of thumb, if $\kappa = \mathcal{O}(10^m)$ for some
$m\in\mathbb{N}_0$ and the sum is evaluated using floating point
arithmetic of $n$ decimal digits accuracy, the result can be
expected to have about $n-m$ correct decimal digits. For the integral
$I$ he reported $\kappa\approx 4.6\,\times\,10^{8}$. He used
standard double precision numbers ($n\approx 16$) for his computation
to obtain:
\begin{equation}
I\approx 181.4393098137807101011276,
\end{equation}
and concluded that $I\approx181.43931$ should be a correct
rounding.

We repeated his experiment using floating point arithmetic of
100 decimal places and exact rational arithmetic for the polynomials.
Due to the different formulation, we obtain a slightly smaller condition number
of $\kappa\approx 1.1\,\times\,10^{8}$. The result should thus be accurate to
about 92 decimal places. The first 50 of those are:
\begin{equation}
I \approx 181.43931117544219248665837073310890818752885155281,
\end{equation} 
thereby precisely confirming Hackbusch's calculations and precision
estimates.

The benefit of the new formulation alone can be seen when setting both
\lstinline{float_type} and \lstinline{number_type} to \lstinline{double}
and using Rump's summation algorithm\autocite{rump2009} for the involved sums.
We then instead obtain $I\approx 181.43931120\dotsc$, giving two more accurate
digits compared to Hackbusch's result. This example is implemented in the
file \lstinline{hackbusch_example.cpp}.

\section{Application to General Data $f\in L^2(\mathbb{R}^3)$}
As an exemplary application, we briefly outline how the results of this work can
be combined with a fast multipole method to obtain an efficient solver for the
Poisson problem $-\Delta u = f$ on the whole-space $\mathbb{R}^3$. We will assume
that $f\in L^2(\mathbb{R}^3)$ is compactly supported. For the sake of both
simplicity and brevity, we only consider uniform discretisations, however
extensions to adaptive schemes are certainly possible.

\subsection{Cartesian Grid and Approximation Spaces}
For arbitrary $\vv{i}=(i_1,i_2,i_3)^\top\in\mathbb{Z}^3$ we define associated
cubes $Q_{\vv{i}}$ of a Cartesian grid $\Omega_h$ of mesh-size $h>0$:
\begin{align}
Q_{\vv{i}} &\coloneqq [hi_1,h(i_1+1)]\times [hi_2,h(i_2+1)]\times [hi_3,h(i_3+1)] \qquad \forall\vv{i}\in\mathbb{Z}^3,\\
\Omega_h   &\coloneqq \bigcup\left\lbrace Q_{\vv{i}} \,\left|\, Q_{\vv{j}}\cap\supp f \neq\emptyset\ 
\forall{\vv{j}}\in\mathbb{Z}^3 \text{ with } |\vv{i}-\vv{j}|_{\infty} \leq 1\right.\right\rbrace.
\end{align}
The definition of $\Omega_h$ ensures that all of its cells $Q_{\vv{i}}$ either
intersect the support of $f$, or have a direct neighbour $Q_{\vv{j}}$ which does.
The multipole method will consider near and far regions for each cell
$Q_{\vv{i}}\subset\Omega_h$:
\begin{align}
\text{Near}(Q_{\vv{i}}) &\coloneqq \bigcup\left\lbrace Q_{\vv{j}}\subset\Omega_h \,\left|\,
|\vv{i}-\vv{j}|_{\infty} \leq 1\right.\right\rbrace, \\
\text{Far}(Q_{\vv{i}})  &\coloneqq \mathbb{R}^3\setminus\text{Near}(Q_{\vv{i}}).
\end{align}

On the mesh $\Omega_h$ we employ spaces of piece-wise polynomials:
\begin{equation}
V_h^n(\Omega_h) \coloneqq \left\lbrace v_h\in L^2(\Omega_h) \,\left|\,
v_h|_{Q_{\vv{i}}}\in\mathbb{P}_n\ \forall Q_{\vv{i}}\subset\Omega_h\right.\right\rbrace,\qquad n\in\mathbb{N}
\end{equation}
where $\mathbb{P}_n$ is the space of all polynomials of \emph{total degree $n-1$
or less}. Because we do not enforce continuity between the individual cells,
orthonormal Legendre polynomials $b_k$, $k\in\lbrace 1,\dotsc,
K\coloneqq\dim V_h^n(\Omega_h)\rbrace$ can be used as a basis for these spaces.
The following definition of index sets then is useful:
\begin{equation}
\text{idx}(Q_{\vv{i}}) \coloneqq \lbrace k\,|\,\supp b_k = Q_{\vv{i}}\rbrace.
\end{equation}

\subsection{Polynomial Approximation of $f$}
It is then straight forward to compute the $L^2(\Omega_h)$-projection $f_h$
of $f$. For $n\geq 2$ the following error bound is standard:\autocite{brenner2008}
\begin{align}
\Vert f-f_h\Vert_{H^{-2}(\mathbb{R}^3)} &\leq
C(n)h^{s+2}\Vert f\Vert_{H^s(\Omega_h)}\qquad 0\leq s\leq n,\\
f_h &\coloneqq \sum_{k=1}^{K} \underbrace{\left( \int_{\mathbb{R}^3}b_k(\vv{x})f(\vv{x})
\,{\mathrm d}\vv{x}\right)}_{\eqqcolon c_k} b_k,
\end{align}
where we implicitly extended $f_h$ by zero outside $\Omega_h$ and used that
$f_h\equiv 0$ on the boundary cells of $\Omega_h$. By the continuity of the
Newton potential $G\star(\bullet):H^{-2}\to L^2$, one then immediately obtains
the following error bound for the resulting potential $u_h\coloneqq G\star f_h$:%
\autocite[Theorem~3.1.2]{sauter2011}
\begin{multline}\label{eqn:error}
\Vert u - u_h\Vert_{L^2(\Omega_h)} = \\
\Vert G\star(f - f_h)\Vert_{L^2(\Omega_h)}   \leq
 C_G \Vert f - f_h\Vert_{H^{-2}(\Omega_h)} \\\leq
 C_GC(n)h^{s+2}\Vert f\Vert_{H^s(\Omega_h)}\qquad 0\leq s\leq n.
\end{multline}
Similar error estimates also hold outside $\Omega_h$, however $u$ and $u_h$
usually do not decay fast enough to zero at infinity to lie in
$L^2(\mathbb{R}^3)$.

\subsection{Polynomial Approximation of $u_h$}
We will see that for $\vv{x}\notin\Omega_h$ the potential $u_h(\vv{x})$ can be
evaluated efficiently using only multipole expansions. On the other hand, for
$\vv{x}\in\Omega_h$, we proceed similarly to Hackbusch\autocite{hackbusch2008}
and compute the $L^2(\Omega_h)$-projection $\tilde{u}_h$ of $u_h$ onto
$V_h^{n+2}(\Omega_h)$, i.\,e., with $L\coloneqq\dim V_h^{n+2}(\Omega_h)$ we compute:
\begin{equation}
\tilde{u}_h = \sum_{l=1}^{L}\underbrace{\left(\int_{\Omega_h}u_h(\vv{x})b_l(\vv{x})
\,{\mathrm d}\vv{x}\right)}_{\eqqcolon d_l} b_l.
\end{equation}
Using standard arguments one obtains that -- up to a constant factor -- the error
bound~\eqref{eqn:error} also holds with $u_h$ replaced by $\tilde{u}_h$. The main
challenge thus lies in the computation of the coefficients $d_l$.

\subsection{Multipole Expansions}
In order to evaluate $u_h$ efficiently, we combine the results of this work with a fast
multipole method. It is well-known that $G$ can be expanded using solid harmonics:%
\autocite[Equation~(46)]{dehnen2014}
\begin{equation}\label{eqn:Gexpansion}
4\pi G(\vv{x}-\vv{y}) \approx
\sum_{p=0}^{P-1}\sum_{q=-p}^{p}\overline{\Theta_p^q(\vv{x})}\Upsilon_p^q(\vv{y}) 
\qquad |\vv{x}|>|\vv{y}|,
\end{equation}
where $P$ is the order of the expansion, which is known to converge rapidly. For
simplicity, we will assume that the error introduced by this approximation is negligible.
In practice this is assumption often is justified for large enough values of $P$, and the
discretisation error~\eqref{eqn:error} is the limiting factor. Rigorous error estimates
for fast multipole methods can be found in the given references, however they are
very technical and beyond the scope of this work.

For each cell $Q_{\vv{i}}\subset\Omega_h$ we may thus define \emph{multipole
moments} $M_{p}^{q}(Q_{\vv{i}})$ as:
\begin{equation}\label{eqn:multipole}
M_{p}^{q}(Q_{\vv{i}}) \coloneqq \frac{1}{4\pi}
\int_{Q_{\vv{i}}}\Upsilon_p^q(\vv{y}-\vv{y}_{\vv{i}})f_h(\vv{y})\,
{\mathrm d}\vv{y}, \qquad \vv{y}_{\vv{i}} \coloneqq \text{centre}(Q_{\vv{i}}).
\end{equation}
Writing $\vv{x}-\vv{y} =(\vv{x}-\vv{y}_{\vv{i}})-(\vv{y}-\vv{y}_{\vv{i}})$,
together with~\eqref{eqn:Gexpansion} we thus immediately obtain:
\begin{equation}\label{eqn:multipole-expansion}
\int_{Q_{\vv{i}}}G(\vv{x}-\vv{y})f_h(\vv{y})\,{\mathrm d}\vv{y}
\approx
\sum_{p=0}^{P-1}\sum_{q=-p}^{p}M_p^q(Q_{\vv{i}})\overline{\Theta_p^q(\vv{x}-\vv{y}_{\vv{i}})}
\qquad\vv{x}\in\text{Far}(Q_{\vv{i}}).
\end{equation}

Note that these moments can be computed efficiently and exactly using standard
quadrature methods, as both $\Upsilon_p^q$ and $f_h$ are polynomials on
$Q_{\vv{i}}$. Even more, writing $f_h = \sum_k c_k b_k$, the mapping
$(c_k)_{k\in\text{idx}(Q_{\vv{i}})} \mapsto (M_{p}^q(Q_{\vv{i}}))_{p,q}$ is linear,
and can thus be stored as a matrix. It is now important to realise that this matrix
is independent of $\vv{i}$: it may thus be \emph{pre-computed} on a single cell, say
$\vv{i}=\vv{0}$, and then reused for all others.

Now remember that $f_h$ vanishes on the cells at $\Omega_h$'s boundary. The
multipole expansions are therefore all we need to evaluate $u_h(\vv{x})$ for
$\vv{x}\notin\Omega_h$. The summation over the individual cells $Q_{\vv{i}}$ and
their respective multipole expansions can be carried out efficiently using so-called
tree-codes.\autocite{barnes1986} The remaining parts of this section therefore
will focus on the case $\vv{x}\in\Omega_h$.

\subsection{Local Expansions}
The fast multipole method\autocite{greengard1987,dehnen2002,dehnen2014} is an
algorithm that takes the multipole moments $M_{p}^{q}(Q_{\vv{i}})$ of all cells
as input and efficiently computes \emph{local moments} $L_{p}^{q}(Q_{\vv{i}})$
for all cells. The details of this algorithm are beyond the scope of this work.
The key property of the local moments is that they locally give the far-field
part with high accuracy:\autocite[Equation~(3a)]{dehnen2014}
\begin{equation}\label{eqn:local}
\int_{\text{Far} (Q_{\vv{i}})}G(\vv{x}-\vv{y})f_h(\vv{y})\,{\mathrm d}\vv{y} \approx
\sum_{p=0}^{P-1}\sum_{q=-p}^{p}L_{p}^{q}(Q_{\vv{i}})\overline{\Upsilon_p^q(\vv{y}_{\vv{i}}-\vv{x})}
\qquad \vv{x}\in Q_\vv{i}.
\end{equation}
This means that the influence of \emph{all cells} $Q_{\vv{j}}\subset
\text{Far}(Q_{\vv{i}})$ is combined in a \emph{single} local expansion and can
thus be evaluated efficiently.

\subsection{Computation of the Coefficients $d_l$}
For each cell $Q_{\vv{i}}\subset\Omega_h$ and each $l\in\text{idx}(Q_{\vv{i}})$
we need to compute:
\begin{multline}
d_l = \int_{Q_{\vv{i}}} u_h(\vv{x})b_l(\vv{x})\,{\mathrm d}\vv{x}
= \int_{Q_{\vv{i}}}\int_{\Omega_h}G(\vv{x}-\vv{y})f_h(\vv{y})b_l(\vv{x})\,{\mathrm d}\vv{y}\,{\mathrm d}\vv{x} = \\
\int_{Q_{\vv{i}}}\int_{\text{Near}(Q_{\vv{i}})}G(\vv{x}-\vv{y})f_h(\vv{y})b_l(\vv{x}) \,{\mathrm d}\vv{y}\,{\mathrm d}\vv{x}
+
\int_{Q_{\vv{i}}}\int_{\text{Far} (Q_{\vv{i}})}G(\vv{x}-\vv{y})f_h(\vv{y})b_l(\vv{x}) \,{\mathrm d}\vv{y}\,{\mathrm d}\vv{x} \\
\eqqcolon  d_{\text{Near},l} + d_{\text{Far},l}.
\end{multline}

For $d_{\text{Far},l}$ we use the local expansion~\eqref{eqn:local}.
Just as in the case of the multipoles $M_p^q(Q_{\vv{i}})$, the mapping
$\bigl(L_p^q(Q_{\vv{i}})\bigr)_{p,q}\mapsto (d_{\text{Far},l})_{l\in\text{idx}(Q_{\vv{i}})}$
can be stored as a matrix that is independent of $\vv{i}$ and can thus be
pre-computed.

For the near-field we write:
\begin{equation}
d_{\text{Near},l} = \sum_{Q_{\vv{j}}\subset\text{Near}(Q_{\vv{i}})} 
\int_{Q_{\vv{i}}}\int_{Q_{\vv{j}}}G(\vv{x}-\vv{y})f_h(\vv{y})b_l(\vv{x}) \,{\mathrm d}\vv{y}\,{\mathrm d}\vv{x} \eqqcolon
\sum_{Q_{\vv{j}}\subset\text{Near}(Q_{\vv{i}})}  d_{Q_{\vv{j}},l}, 
\end{equation}
and note that there are at most $3^3 = 27$ cells $Q_{\vv{j}}$ in this sum. Let us
now fix one such $\vv{j}$. We have with $f_h = \sum_k c_k b_k$:
\begin{equation}
d_{Q_{\vv{j}},l} = \sum_{k\in\text{idx}(Q_{\vv{j}})}c_k
\int_{Q_{\vv{i}}}\int_{Q_{\vv{j}}}G(\vv{x}-\vv{y})b_k(\vv{y})b_l(\vv{x}) \,{\mathrm d}\vv{y}\,{\mathrm d}\vv{x}
\end{equation}
This double integral can be evaluated \emph{exactly} using the methods described in
this paper. However, this evaluation is costly. It is therefore important to notice that
the mapping $(c_{k})_{k\in\text{idx}(Q_{\vv{j}})}\mapsto
(d_{Q_{\vv{j}},l})_{l\in\text{idx}(Q_{\vv{i}})}$ can again be stored as a matrix,
and this matrix only depends on the \emph{difference} $\vv{j}-\vv{i}$. We can thus
pre-compute these so-called \emph{interaction matrices} for $\vv{i}=\vv{0}$ and the
corresponding 27 values of $\vv{j}$. Moreover, changing the mesh-size $h$ only
changes these matrices by a constant factor. It is thus sufficient to pre-compute
these matrices once for $h=1$, using very high precision to avoid cancellation, and
then store the result in, e.\,g., double precision.

Thus, after the pre-computation of all involved matrices, and after the
coefficients $c_k$ of $f_h$ have been computed, the entire algorithm reduces to
a set of small to moderately sized, dense matrix-vector multiplications which
can be efficiently carried out using the BLAS.

\subsection{Numerical Example}
To illustrate the efficiency of the above approach, we prescribe the solution
$u$ of $-\Delta u = f$ to be Friedrichs's mollifier:
\begin{equation}
u: \mathbb{R}^3\to\mathbb{R},\qquad \vv{x}\mapsto
\begin{cases}
\exp\left(-\frac{1}{1-|\vv{x}|^2}\right) & \text{ if } |\vv{x}| < 1, \\
0                                        & \text{ else. }
\end{cases}
\end{equation}

We conduct a simple convergence study using $n=4$ for the approximations
$f_h\in V_h^n(\Omega_h)$ of $f$ and $\tilde{u}_h\in V_h^{n+2}(\Omega_h)$ of $u$,
and vary $h$. For the fast multipole method we use a simple code employing
a fixed expansion of order $P = 20$, without any acceleration of the so-called
M2L-operator. A code with variable, adaptive expansion orders and M2L
accelerations would certainly result in further speed-up,\autocite{dehnen2014}
but lies beyond the scope of this work.

For the computations we use a simple, low-budget laptop with an Intel Core
i5-7200U processor. The code is available separately from the library, under
the terms and conditions of the (full) GNU General Public License, version
three or later, as supplementary material to this paper.

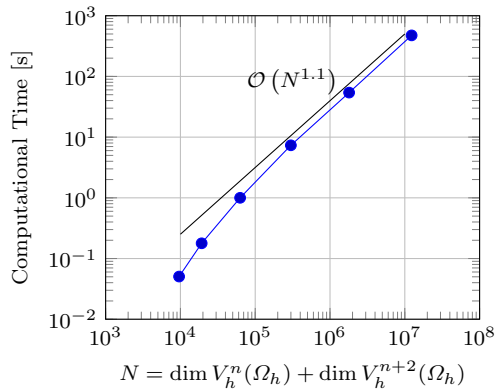
\begin{figure}
\centering
\begin{tikzpicture}
\begin{loglogaxis}
[
	small,
    grid = major,
	xmin = {1e3},
	xmax = {1e8},
    ymin = {1e-2},
	ymax = {1e3},
    xlabel = {$N =\dim V_h^n(\Omega_h) + \dim V_h^{n+2}(\Omega_h)$},
    ylabel = {Computational Time [s]},
    legend pos = south west,
]
\addplot table[x=NDOF,y=time]{timings.dat};
\addplot[domain=1e4:1e7]{1e-5*x^1.1} node[above=12.0,pos=.5] {$\bigO{N^{1.1}}$};
\end{loglogaxis}
\end{tikzpicture}
\caption{\label{fig:timings}The computational time scales almost linearly
in the number of unknowns. Note that this code was executed on a low-budget
laptop computer.}
\end{figure}

\begin{figure}
\centering
\begin{tikzpicture}
\begin{semilogyaxis}
[
	small,
    grid = major,
	xmin = {0},
	xmax = {5},
    ymin = {1e-8},
	ymax = {1e1},
    ylabel = {$L^2(\Omega_h)$-Error},
    xlabel = {$h=2^{-k}$},
    legend pos = south west,
    ytick = { 1e-8, 1e-7, 1e-6, 1e-5, 1e-4, 1e-3, 1e-2, 1e-1, 1e-0, 1e1 },
    minor ytick = { 2e-8, 4e-8, 6e-8, 8e-8,
                    2e-7, 4e-7, 6e-7, 8e-7,
                    2e-6, 4e-6, 6e-6, 8e-6,
                    2e-5, 4e-5, 6e-5, 8e-5,
                    2e-4, 4e-4, 6e-4, 8e-4,
                    2e-3, 4e-3, 6e-3, 8e-3,
                    2e-2, 4e-2, 6e-2, 8e-2,
                    2e-1, 4e-1, 6e-1, 8e-1,
                    2e-0, 4e-0, 6e-0, 8e-0  },
]
\addplot table[x=k,y=fh]{errors.dat};
\addplot table[x=k,y=uh]{errors.dat};
\addplot[domain=3.5:6]{2.5e4*2^(-4*x)} node[above=5.0,pos=0.25] {$\bigO{h^{4}}$};
\addplot[domain=2.5:6]{5e2*2^(-6*x)} node[above=5.0,pos=.4] {$\bigO{h^{6}}$};
\addlegendentry{$\Vert f-f_h\Vert_{L^2(\Omega_h)}$};
\addlegendentry{$\Vert u-\tilde{u}_h\Vert_{L^2(\Omega_h)}$};
\end{semilogyaxis}
\end{tikzpicture}
\caption{\label{fig:errors}$L^2(\Omega_h)$-errors of the approximations
$f_h$ and $u_h$ when solving $-\Delta u = f$, using order $n=4$ as the
grid-size $h$ tends to zero. While for $f_h$ we can only see the beginning
of the asymptotic behaviour, $\tilde{u}_h$ converges rapidly at sixth order
as predicted by the error-bound~\eqref{eqn:error}.}
\end{figure}
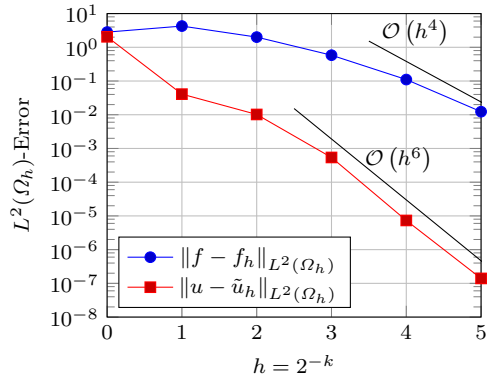

Already this simple code achieves reasonable performance. The fast multipole
method has a theoretical time complexity of $\mathcal{O}(N)$. As can be seen in
Figure~\ref{fig:timings}, the code has an empirical complexity that almost matches
this result. Figure~\ref{fig:errors} shows the errors of $f_h$ and $\tilde{u}_h$
as $h\to 0$. While for $f_h$ we can only observe the beginning of the asymptotic
range, the sixth order convergence for $\tilde{u}_h$ is already clearly visible.

We can thus conclude that the combination of exact integration formulæ with
a fast multipole method allows us to compute the Newton potential of general
data $f$ efficiently and to very high accuracy.

\printbibliography
\end{document}